\documentclass[amscd,amssymb,verbatim,a4paper]{amsart}

\theoremstyle{plain}
\newtheorem{theorem}{Theorem}[section]
\newtheorem{lemma}[theorem]{Lemma}

\theoremstyle{definition}

\theoremstyle{remark}
\newtheorem{remark}[theorem]{Remark}
\newtheorem{corollary}[theorem]{Corollary}

\newtheorem{notations}[theorem]{Notations}

\newcommand{\ic}{\ensuremath{\mathcal{I}}}
\newcommand{\oc}{\ensuremath{\mathcal{O}}}

\newcommand{\lc}{\ensuremath{\mathcal{L}}}

\newcommand{\nc}{\ensuremath{\mathcal{N}}}

\newcommand{\mc}{\ensuremath{\mathcal{M}}}

\newcommand{\Pq}{\mathbb{P}^4}
\newcommand{\Pcq}{\mathbb{P}^5}
\newcommand{\Psx}{\mathbb{P}^6}

\newcommand{\Pn}{\mathbb{P}^n}
\newcommand{\bZ}{\mathbb{Z}}

\begin{document}

\title[On smooth divisors of a projective hypersurface.]{On smooth divisors of a projective hypersurface.}

\author{Ellia Ph.}
\address{Dipartimento di Matematica, via Machiavelli 35, 44100 Ferrara (Italy)}
\email{phe@dns.unife.it}

\author{Franco D.}
\address{Dipartimento di Matematica e Applicazioni "R. Caccioppoli", Univ. Napoli "Federico II", Ple Tecchio 80, 80125 Napoli (Italy)}
\email{davide.franco@unina.it}

\date{16/06/2004}

\maketitle

\hskip7cm{\it Dedicated to Christian Peskine.}
\bigskip

\section*{Introduction.}

This paper deals with the existence of smooth divisors of a projective hypersurface $\Sigma \subset\Pn $ (projective space over an algebraically closed field of characteristic zero). According to a celebrated conjecture of Hartshorne, at least when $n\geq 7$, any such a variety should be a complete intersection.
Since the existence of smooth, non complete intersection, subcanonical $X \subset \Pn$ of codimension two
is equivalent, via the correspondance of Serre, to the existence of indecomposable rank two vector bundles on $\Pn$ and since no indecomposable vector bundle of $\Pn $, $n\geq 5$, is presently known,
it is widely believed that any smooth, subcanonical subvariety of $\Pn $, $n\ge5$, of codimension two is a complete intersection. Furthermore recall that, by a theorem of Barth, the subcanonical condition is automatically satisfied if $n \geq 6$. This in turn implies that a smooth (subcanonical if $n=5$) divisor of a
projective hypersurface $\Sigma \subset\Pn $, $n\geq 5$, is a complete intersection too.
\par In this paper we show that, roughly speaking, for any $\Sigma \subset \Pn$ there can be at most $\textit{finitely many}$ exceptions to the last statement. Indeed our main result is:

\begin{theorem}
\label{mainthm}
Let $\Sigma \subset \Pn$, $n \geq 5$ be an integral hypersurface of degree $s$. Let $X \subset \Sigma$ be a smooth variety with $dim(X)=n-2$. If $n=5$, assume $X$ subcanonical. If $X$ is not a complete intersection in $\Pn$, then:
$$d(X) \leq \frac{s(s-1)[(s-1)^2-n+1]}{n-1}+1.$$
\end{theorem}

In other words a smooth codimension two subvariety of $\Pn$, $n \geq 5$ (if $n=5$, we assume $X$ subcanonical) which is not a complete intersection cannot lie on a hypersurface of too low degree (too low with respect to its own degree) and, {\it on a fixed hypersurface}, Hartshorne's conjecture in codimension two is "asymptotically" true.
\par
The starting point is Severi-Lefschetz theorem which states that if $n \geq 4$ and if $X$ is a Cartier divisor on $\Sigma$, then $X$ is the complete intersection of $\Sigma$ with another hypersurface. For instance if $\Sigma $ is either smooth or singular in a finite set of points and if $n \geq 5$, the picture is very clear:
\begin{enumerate}
\item
there $\textit{exists}$ smooth $X\subset \Sigma$ with $dim(X)=n-2$ and with degree arbitrarily large;
\item
$\textit{any}$ smooth $X\subset \Sigma$ with $dim(X)=n-2$ $\textit{is}$ a complete intersection of $\Sigma$ with another hypersurface
\item $\textit{no}$ smooth $X\subset \Sigma$ with $dim(X)=n-2$ can meet the singular locus of $\Sigma$.
\end{enumerate}

\vskip1cm

Using Theorem \ref{mainthm} we get (the first statement comes again from an easy application of the Theorem of Severi-Lefschetz-Grothendieck):

\begin{theorem}
\label{sigma}  
Let $\Sigma \subset \Pn$, $n\geq 5$, be an integral hypersurface of degree $s$
with  \par $dimSing(\Sigma)\geq1$.
\begin{enumerate}
\item If $n\geq 6$ and
$dimSing(\Sigma)\leq n-5$ then $\Sigma $ does not contain any smooth variety of dimension $n-2$.
\item Suppose $dimSing(\Sigma)\geq n-4$. If $X\subset \Sigma $ is smooth, subcanonical, with $dim(X)=n-2$ then
$d(X)\leq s\frac{(s-1)((s-1)^2-n+1)}{n-1}+1$.
\end{enumerate}
\end{theorem}

We point out a consequence of this result.

\begin{corollary}
\label{sigmahilb} Let $\Sigma \subset \Pn$, $n\geq 5$, be an integral hypersurface
s.t. $dimSing(\Sigma)\geq 1$.
\begin{enumerate}
\item If $n\geq 6$ and
$dimSing(\Sigma)\leq n-5$ then $\Sigma $ does not contain any smooth variety of dimension $n-2$.
\item Suppose $dimSing(\Sigma)\geq n-4$. Then there are only finitely many components of
$\mathcal{H}ilb(\Sigma)$ containing  smooth, subcanonical varieties of dimension $n-2$.
\end{enumerate}
\end{corollary}

Last but not least, at the end of the paper we show how this circle of ideas  allows to improve the main results of \cite{EF} about subcanonical varieties of $\Pcq $ and $\Psx$:

\begin{theorem}
\label{n=s=5}
Let $X \subset \Pcq$ be a smooth threefold with $\omega _X \simeq \oc _X(e)$. If $h^0(\ic _X(5)) \neq 0$,
then $X$ is a complete intersection.  
\end{theorem}

\begin{theorem}
\label{n=s=6}
Let $X \subset \Psx$ be a smooth fourfold. If $h^0(\ic _X(6)) \neq 0$,
then $X$ is a complete intersection.  
\end{theorem}

Theorem \ref{mainthm} follows, thanks to a crucial remark essentially proved in \cite{EP} (see Lemma \ref{l1}), from a bound of $e$ (where $\omega _X \simeq \oc _X(e)$), see Theorem \ref{thmSpec}, which can be viewed as a strong (since the degree is not involved) generalization of the "Speciality theorem" of Gruson-Peskine \cite{GP}. The proof of this bound is quite simple if $X \cap Sing(\Sigma )$ has the right dimension. This is done in the first section where a weaker version of Theorem \ref{thmSpec} and hence of Theorem \ref{mainthm} is proved (if $n=5$ we assume $Pic(X) \simeq \bZ .H$). In the second section we show how a refinement of the proof yields our final result. Finally let's observe that our approach doesn't apply to the case $n=4$.

\vskip1cm

\textbf{Acknowledgment:} It is a pleasure to thank Enzo Di Gennaro who explained to one of us (D.F.) some of the deep results  of \cite{K}. 

\section{Reduction and the speciality theorem, weak version.}

\begin{notations}
Given a projective scheme $Y \subset \Pn $ we denote by $d(Y)$ the \textit{degree} of $Y$.
\end{notations}

\begin{notations}
\label{not2} In this section, $X \subset \Pn, n \geq 5$, will denote a smooth, non degenerate, codimension two
subvariety which is not a complete intersection. We will always assume $X$ subcanonical: $\omega _X \simeq
\oc _X(e)$; notice that this condition is fullfilled if $Pic(X) \simeq \bZ .H$; finally, thanks to a theorem of Barth, this last condition is automatically fullfilled if $n \geq 6$.\par\noindent
By Serre's construction we may associate to $X$ a rank two vector bundle:
$$ 0 \to \oc \to E \to \ic _X(e+n+1) \to 0 $$
The Chern classes of $E$ are: $c_1(E)=e+n+1,c_2(E)=d(X)=:d$.\par\noindent
Let $\Sigma$ be an hypersurface of degree $s $ containing $X$.
Then $\Sigma $ gives a section of $\ic _X(s)$ which lifts to a section  $\sigma _{\Sigma}\in H^0(E(-e-n-1+s))$ (notice that $\sigma_{\Sigma }$ is uniquely defined if $e+n+1-s<0$).
{\it Assume} that $Z$, the zero-locus of $\sigma_{\Sigma }$, has codimension two. Notice that since $X$ is not a complete intersection, this certainly holds if $s = min\{t\:|\:h^0\ic _X(t)) \neq 0\}$.
Anyway, if $Z$ has codimension two, then $d(Z)=c_2(E(-e-n-1+s))=d-s(e+n+1-s)$ and $\omega _Z \simeq \oc _Z(-e-2n-2+2s)$. 
\end{notations}

\begin{remark}
\label{snb} By \cite{R}, if $X\subset \Sigma \subset \Pn $, $n\geq 3$, with $\omega _X \simeq \oc _X(e)$ and $d(\Sigma)\leq n-2$ then $X$ is complete intersection, hence in the remainder of this paper we will assume $s\geq n-1$.
\end{remark}

\begin{remark}
\label{jac} Notice that $E(-e-n-1)\mid _X\simeq \nc ^*_X$. It is well known 
that the scheme $X\cap Z$ is the base locus of the jacobian system of $\Sigma $ on $X$: $X\cap Z=X\cap Jac(\Sigma)$.
So, the \textit{fundamental cycle} (\cite{F} 1.5) of $Z$ in $\mathcal{A}_*(X)$ is $c_2(\nc ^*_X(s))$ as soon as $X$ and $Z$ intersect in the expected codimension.
\end{remark}

The main goal of this section is to prove:

\begin{theorem}[Speciality theorem, weak version]
\label{thmSpecW}
Let $X \subset \Pn$, $n \geq 5$ be a smooth codimension two subvariety. If $n=5$ assume $Pic(X) \simeq \bZ .H$. Let $\Sigma$ be an hypersurface of degree $s$ containing $X$. If $X$ is not a complete intersection, then:
$$e \leq \frac{(s-1)[(s-1)^2-n+1]}{n-1}-n+1$$
where $\omega _X \simeq \oc _X(e)$.
\end{theorem}

Let's see how this is related with a bound of the degree. First recall the following:

\begin{lemma}
\label{l1}
Let $X \subset \Pn$, $n \geq 4$, be a smooth codimension two subvariety which is not a complete intersection. Let $\Sigma$ be an hypersurface of minimal degree containing $X$. Set $s:=d(\Sigma )$.
\begin{enumerate}
\item $n-4 \leq dim(X \cap Sing(\Sigma )) \leq n-3$.
\item If $\omega _X \simeq \oc _X(e)$, then $d(X) \leq s(n-1+e)+1$.
\item If $dim(X \cap Sing(\Sigma ))=n-3$ and if $Pic(X) \simeq \bZ .H$, then $d(X) \leq (s-2)(n-1+e)+1$.
\end{enumerate}
\end{lemma}

\begin{proof}
The first item is \cite{EF}, Lemma 2.1; 2) is \cite{EF} Lemma 2.2 (i) and the last item is \cite{EF} Lemma 2.2 (ii) with $l=2$ (thanks to Severi and Zak theorems $h^1(\ic _X(1))=0$, \cite{Z}). 
\end{proof}

Theorem \ref{thmSpecW} and the second item of this lemma give us immediately:

\begin{theorem}
\label{thmA}
Let $\Sigma \subset \Pn$, $n \geq 5$, be an integral hypersurface of degree $s$. Let $X \subset \Sigma$ be a smooth subvariety with $dim(X)=n-2$. If $n=5$ assume $Pic(X) \simeq \bZ .H$. If $X$ is not a complete intersection, then $d(X) < \frac{s(s-1)[(s-1)^2-n+1]}{n-1}+1$.
\end{theorem}

In order to prove Theorem \ref{thmSpecW} we need some preliminary results.

\begin{lemma}
\label{Ysubcanonical}
Let $\Sigma$ denote an hypersurface of degree $s$ containing $X$. 
With assumptions ($codim(\sigma _{\Sigma})_0=2$) and notations as in \ref{not2}, assume $dim(X\cap Z)=n-4$. Then $Y:=X\cap Z$ is a subcanonical, l.c.i. scheme with $\omega_Y\simeq\oc_Y(2s-n-1)$. Moreover $Y$ is the base locus of the jacobian system of $\Sigma$ in $X$.
\end{lemma}

\begin{proof} 
We are assuming that $Y$ is a proper intersection between $X$ and $Z$ hence
$$
0 \to \oc \to E\mid_X(-e-n-1+s) \to \ic _{Y, X}(-e-n-1+2s) \to 0
$$
so $\nc^*_{Y,X}\simeq E\mid_X(-s)$ and the first statement follows by adjunction. For the last statement, use \ref{jac}.
\end{proof}

\begin{notations} Keep the assumptions of Lemma \ref{Ysubcanonical} and denote by $\Sigma_1$ and $\Sigma_2$ two general partials of $\Sigma$. Since $dim(X\cap Z)=n-4$, $C:= X\cap \Sigma_1\cap \Sigma_2$ is a subcanonical, l.c.i. scheme  containing $Y$ such that $\nc_{C,X}\simeq \oc _X(s-1) \oplus \oc_X(s-1)$. We have $\omega_C\simeq \oc_C(e+2s-2)$. The scheme $C$ is a complete intersection in $X$ which links $Y$ to another subscheme.\\
\end{notations}

\begin{lemma}
\label{R}
With notations as in Lemma \ref{Ysubcanonical}, denote by $R$ the residual to $Y$ with respect to $C$. Then $C=Y\cup R$ is a geometric linkage and $\Delta:= R\cap Y$ is a Cartier divisor of $Y$ such that: $\ic _{\Delta , Y}\simeq \oc_Y (-e-n+1)$.\\
Furthermore: $d(\Delta)\leq (s-1)d(X)((s-1)^2-d(Z))$ and:\\
 $d(Z)(e+n+1) \leq (s-1)[(s-1)^2-d(Z)]$.

\end{lemma}

\begin{proof}
Denote by $Y_{red}$ the support of $Y$ and set
$Y_{red}=Y_1 \cup \dots \cup Y_r$ where $Y_i$, $1\leq i \leq r$, are the irreducible components of $Y_{red}$.
Furthermore, denote by $P_i$ the general point of $Y_i$. Since $Y$ is l.c.i. in $X$ and since $\ic _{Y, X}(s-1)$
is globally generated by the partials of $\Sigma$,  we can find two general elements in $Jac(\Sigma)$ generating the fibers of
$\nc^* _{Y, X}(s-1)$ at each $P_i$, $1\leq i\leq r$. This implies that $R\cup Y$ is a geometric linkage.
\par Now consider the local Noether sequence (exact sequence of liaison):
$$
0\to \ic_C \to \ic _R \to \omega_Y \otimes \omega_C^{-1}\to 0.
$$
we get $$
\omega_Y \otimes \omega_C^{-1}\simeq\frac{\ic_R}{\ic_C}\simeq \frac{\ic_R +\ic_Y}{\ic_C+\ic _Y}\simeq
\frac{\ic_{\Delta}}{\ic_Y}\simeq \ic _{\Delta , Y}$$ (the second isomorphism follow by geometric linkage, since $\ic_R\cap\ic_Y =\ic_C$)
hence
$\omega_Y \otimes \omega_C^{-1}\simeq \oc_Y(-e-n+1)\simeq \ic _{\Delta , Y}$ and we are done.\\
For the last statement, the scheme $\Delta \subset R $ is the base locus of the jacobian system of
$\Sigma $ in $R$, hence $\Delta \subset \tilde{\Sigma}\cap R $ with $\tilde{\Sigma }$ a general element of $Jac(\Sigma)$
and $d(\Delta )\leq d(R)\cdot (s-1)$. We conclude since $d(R)\cdot (s-1)=(d(C)- d(Z))\cdot (s-1)=
((s-1)^2d(X)-d(Z)d(X))\cdot (s-1)$. The last inequality follows from $d(\Delta )=d(Y)\cdot (e+n+1)=d(X)\cdot d(Z) \cdot (e+n+1)$.
\end{proof}

Now we can conclude the proof of Theorem \ref{thmSpecW} (and hence of Theorem \ref{thmA}).

\begin{proof}[Proof of Theorem \ref{thmSpecW}]
It is enough to prove the theorem for $s$ minimal. Let $\Sigma$ be an hypersurface of minimal degree containing $X$, we set $s:=d(\Sigma )$ and $d:=d(X)$. According to Lemma \ref{l1} we distinguish two cases.\\
1) $dim(X \cap Sing(\Sigma ))=n-3$. In this case, by Lemma \ref{l1}, we have $d \leq (s-2)(n-1+e)+1$. On the other hand $d(Z) = d-s(e+n+1-s)$ (see \ref{not2}). It follows that: $d(Z) \leq (s-1)^2 -2(n-1+e)$. Since $d(Z) \geq n-1$ by \cite{R}, we get: $\frac{(s-1)^2-n+1}{2}-n+1 \geq e$. One checks (using $s \geq n-1$) that this implies the bound of Theorem \ref{thmSpecW}.\\
2) $dim(X \cap Sing(\Sigma ))=n-4$. By the last inequality of Lemma \ref{R}, $e \leq (s-1)[\frac{(s-1)^2}{d(Z)}-1]-n+1$. Since $d(Z) \geq n-1$ by \cite{R}, we get the result. 
\end{proof}

\section{The speciality theorem.}

In this section we will refine the proof of Theorem \ref{thmSpecW} for $n=5$ in order to prove Theorem \ref{mainthm} of the introduction. For this we have to assume only that $X$ is subcanonical, which, of course, is weaker than assuming $Pic(X) \simeq \bZ .H$. The assumption $Pic(X) \simeq \bZ . H$ is used just to apply the last statement of Lemma \ref{l1} in order to settle the case $dim(X \cap Sing(\Sigma ))=n-3$. Here instead we will argue like in the proof of the case $dim(X \cap Sing(\Sigma ))=n-4$, but working modulo the divisorial part (in $X$) of $X \cap Sing(\Sigma )$; this will introduce some technical complications, but conceptually, the proof runs as before. Since the proof works for every $n \geq 5$ we will state it in this generality giving thus an alternative proof of Theorem \ref{thmSpecW}.

\begin{notations}
\label{not3}
In this section, with assumptions and notations as in \ref{not2}, we will assume furthermore that $dim(X\cap Z)=n-3$ and  will denote by $L$ the dimension $n-3$ part of $X\cap Z\subset X$; moreover we set $\lc = \oc _X(L) $.
\par\noindent
Set $Y':=res_L(X\cap Z)$, we have $\ic _{Y', X}:=(\ic _{X\cap Z , X}:\ic _{L , X})$. Since we have:
$$0 \to \oc \to E\mid_X(-e-n-1+s)\otimes \lc ^* \to \ic _{Y' , X}(-e-n-1+2s)\otimes (\lc ^*)^2 \to 0$$
it follows that $\nc^*_{Y',X}\simeq E\mid_X(-s)\otimes \lc$ and $Y'$ is a  l.c.i. scheme with
$\omega_{Y'}\simeq\oc_Y(2s-n-1)\otimes (\lc^*)^2$.
\par \noindent
Denote by $\Sigma_1$ and $\Sigma_2$ two general partials of $\Sigma$. Since $X \cap Z = X \cap Sing(\Sigma )$, $\Sigma_1$ and $\Sigma_2$ both contain $L$. Let $C':=res_L(X\cap  \Sigma_1\cap \Sigma_2)$. Since $\nc_{C',X}\simeq (\oc _{C'}(s-1) \oplus \oc_{C'}(s-1))\otimes \lc^*$. We have
$\omega_{C'}\simeq \oc_{C'}(e+2s-2)\otimes (\lc^*)^2$.
\end{notations}

\begin{lemma}
\label{R2}
Denote by $R'$ the residual to $Y'$ with respect to $C'$. Then $C'=Y'\cup R'$ is a geometric linkage and $\Delta':= R'\cap Y'$ is a Cartier divisor of $Y'$ such that: 
$\ic _{\Delta' , Y'}\simeq \oc_{Y'}(-e-n+1)$.
\end{lemma}

\begin{proof}
We argue as in the proof of Lemma \ref{R}: denote by $Y_{red}'$ the support of $Y'$, set
$Y_{red}'=Y_1' \cup \dots \cup Y_r'$, where $Y_i'$, $1\leq i \leq r$, are the irreducible components of $Y_{red}'$, and denote by $P_i$ the general point of $Y_i'$. Choose the partials $\Sigma_1$ and $\Sigma _2$ in such a way that they generate  the ideal sheaf of $X\cap Z$ at each $P_i$, $1\leq i\leq r$. In order to check that $R'\cup Y'$ is a geometric linkage we only need to consider the components contained in $L$. Consider a point $P_i\in L$. Since 
$L\subset X\cap Z \subset \Sigma_1 \cap \Sigma _2$, the local equations of $X\cap Z$  in
$(\ic _{Y, X}(s-1))_{P_i}$ have the form $(lf,lg)$ where $l$ is the equation of $L$,
$lf$ is the equation of $\Sigma_1$ and $lg$  the equation of $\Sigma_2$.
Since $Y':=res_L(X\cap Z)$  and $C':=res_L(X\cap  \Sigma_1\cap \Sigma_2)$ then the ideals of both $Y'$ and $C'$ at $P_i$ are equal to $(f,g)\subset (\ic _{Y, X}(s-1))_{P_i}$.
This implies that $R'\cup Y'$ is a geometric linkage and the remainder of the proof is similar as above.
\end{proof}

\begin{lemma}
\label{lemmaN-3}
Let $\Sigma \subset \Pn$, $n\geq 5$, be an hypersurface of degree $s$ containing $X$, a smooth variety with $dim(X)=n-2$ and $\omega _X \simeq \oc _X(e)$. Assume $\sigma _{\Sigma}$ vanishes in codimension two and $dim(X \cap Sing(\Sigma ))=n-3$ (see \ref{not2}). Then $e < s-n$ or $d(Z)\cdot (e+n+1) \leq (s-1)[(s-1)^2-d(Z)]$.
\end{lemma}

\begin{proof}
We keep back the notations of \ref{not3}. Notice that the fundamental cycle of $Y'$ in $\textbf{A} _{n-4}(X)$ is 
$$c_2(E\mid_X(-e-n-1+s)\otimes\lc^*)=d(Z)H^2 + (e+n+1-2s)H\cap L +L^2\:\:(+)$$ ($H$ represents the hyperplane class and $\cap $ denotes the \textit{cap product} in $\textbf{A}_*(X)$. By abuse of notations, for any  $A\in \textbf{A} _{i}(X)\subset \textbf{A}_*(X)$ we denote by 
$d(A)\in \mathbb{Z}$ the
\textit{degree} of $A$: $d(A):= d(A\cap H^i)$, $A\cap H^i\in A_0(\Pn )\simeq \mathbb{Z}$. \item For any closed subscheme $\Gamma \subset X$ we still denote by $\Gamma \in \textbf{A} _{*}(X)$ the \textit{fundamental cycle}
of $\Gamma $ (\cite{F} 1.5).\\
We claim that: 
$$d(\Delta')\leq (s-1)d(X)((s-1)^2-d(Z))-[(s-1)(e+n-1)+(s-1)^2-d(Z)]d(H^2\cap L)+$$
$$+(e+n-1)d(H\cap L^2)\:\:(*)
$$
Assume the claim for a while and let's show how to conclude the proof.
Combining \ref{R2} with $(*)$ we get
$$
d(\Delta')=d(Y')(e+n-1)\leq
$$
$$
\leq (s-1)d(X)((s-1)^2-d(Z))-[(s-1)(e+n-1)+(s-1)^2-d(Z)]d(H^2\cap L)+$$
$$+(e+n-1)d(H\cap L^2)
$$ and by $(+)$ above 
$$
d(\Delta')=(e+n-1)d(H\cap(d(Z)H^2 + (e+n+1-2s)H\cap L +L^2))\leq
$$
$$
\leq (s-1)d(X)((s-1)^2-d(Z))-[(s-1)(e+n-1)+(s-1)^2-d(Z)]d(H^2\cap L)+$$
$$+(e+n-1)d(H\cap L^2).
$$
If $e<s-n$ we are done, so we can assume $e+n\geq s$.
We have 
$$d(X)d(Z)(e+n-1)\leq (s-1)d(X)((s-1)^2-d(Z))+$$
$$+[(e+n-1)(s-e-n)-(s-1)^2+d(Z)]d(L)$$

To conclude it is enough to check that $(e+n-1)(s-e-n)-(s-1)^2+d(Z)\leq 0$. Since $d(Z)=d-s(e+n+1-s)$ (see \ref{not2}) and since $d \leq s(n-1+e)+1$ by Lemma \ref{l1}, this follows from: $s(n-1+e)+1 \leq s(e+n+1-s)+(s-1)^2+(e+n-s)(e+n-1)$. A short computation shows that this is equivalent to $0 \leq (e+n-s)(e+n-1)$, which holds thanks to our assumption $e+n\geq s$.\\
{\it Proof of the claim:}\\
Denote by $\mid M \mid $ the moving part of the Jacobian system of $\Sigma $ in $X$ and by $\mc $ the corresponding line bundle.
The scheme $\Delta '$ is the base locus of $\mid M \mid_{R'}$ hence 
$\Delta '\subset \tilde{M}\cap R'$ where $\tilde{M}$ is a general element of $\mid M \mid $. We have
$$
d(\Delta ')\leq d(\tilde{M}\cap R')=d(c_1(\mc _{R'})).
$$
\par
In order to prove the statement we need to calculate the cycle 
$c_1(\mc _{R'})\in \textbf{A} _{n-5}(X)$.
First of all we calculate the fundamental cycle of $R'$ in $\textbf{A} _{n-4}(X)$:
$$R'\sim C'-Y'\sim ((s-1)H-L)^2-(d(Z)H^2 + (e+n+1-2s)H\cap L +L^2)=
$$
$$
=((s-1)^2-d(Z))H^2-(e+n-1)H\cap L.$$
Finally, the cycle $c_1(\mc _{R'})\in \textbf{A} _{n-5}(X)$ is:
$$c_1(\mc _{R'})\sim ((s-1)H-L)\cap R'\sim 
$$
$$
\sim (s-1)((s-1)^2-d(Z))H^3-((s-1)(e+n-1)+(s-1)^2-d(Z))H^2\cap L+(e+n-1)H\cap L^2.$$
The claim follows from:
$$
d(\Delta')\leq d(c_1(\mc _{R'}))=
$$
$$
d((s-1)((s-1)^2-d(Z))H^3-((s-1)(e+n-1)+(s-1)^2-d(Z))H^2\cap L+(e+n-1)H\cap L^2)
$$
\end{proof}

Now we can state the improved version of Theorem \ref{thmSpecW}:

\begin{theorem}[Speciality theorem]
\label{thmSpec}
Let $X \subset \Pn$, $n \geq 5$, be a smooth variety with $dim(X)=n-2$ and $\omega _X \simeq \oc _X(e)$. Let $\Sigma \subset \Pn$ denote an hypersurface of degree $s$ containing $X$. If $X$ is not a complete intersection, then:
$$e \leq \frac{(s-1)[(s-1)^2-n+1]}{n-1}-n+1.$$
\end{theorem}

\begin{proof}
It is sufficient to prove the theorem for $s$ minimal. We distinguish two cases (see Lemma \ref{l1}).\\
If $dim(X \cap Sing(\Sigma ))=n-4$, then we argue exactly as in the proof of Theorem \ref{thmSpecW}.\\
If $dim(X \cap Sing(\Sigma ))=n-3$, then by Lemma \ref{lemmaN-3} we have $e < s-n$ or $d(Z)\cdot (e+n+1) \leq (s-1)[(s-1)^2-d(Z)]$. In the first case we conclude using $s \geq n-1$ (Remark \ref{snb}) and, in the second case, we conclude using the fact that $d(Z) \geq n-1$ by \cite{R}.
\end{proof}

\begin{proof}[Proof of Theorem \ref{mainthm}]
As explained in the Section 1, it follows from Theorem \ref{thmSpec} and Lemma \ref{l1}.
\end{proof}

\section{Proofs of \ref{sigma} and of \ref{sigmahilb}.}

\begin{proof}[Proof of Theorem \ref{sigma}] If $X$ is not a complete intersection, this follows from Theorem \ref{mainthm}. Assume $X$ is a complete intersection. Let $F$ and $G$ ($d(F)=f,d(G)=g$) be two generators of the ideal of $X$. Then the equation of $\Sigma$ has the form $PF+QG$. But since $\Sigma$ is irreducible and since $X \cap Sing(\Sigma )\neq \emptyset$, then both $P$ and $Q$ have degree $>0$. This implies $s-1 \geq f$ and $s-1 \geq g$ hence $d=fg \leq (s-1)^2 < s\frac{(s-1)((s-1)^2-n+1}{n-1}+1$.
\end{proof}

\begin{proof}[Proof of Corollary \ref{sigmahilb}] The argument goes as in the proof of  \cite{CDG} Lemma 4.3: by
\cite{K} the coefficients of the Hilbert polynomial of $X$ can be bounded in terms of the degree $d$ hence in terms of $s$, by \ref{sigma}, and there are finitely many components of $\mathcal{H}ilb(\Sigma)$ containing  smooth varieties of dimension $n-2$.
\end{proof}

\section{Proof of \ref{n=s=5} and \ref{n=s=6}}

\begin{notations} By \cite{EF}, we may assume that $X$ lies on an irreducible hypersurface $\Sigma$ of degree $n$, $5 \leq n \leq 6$ and that $h^0(\ic _X(n-1))=0$. The assumption of \ref{not2} is satisfied and by Lemma \ref{R} and Lemma \ref{lemmaN-3}, we get: $e <s-n$ or $d(Z)\cdot (e+n-1) \leq (s-1)[(s-1)^2-d(Z)]$. The first case cannot occur in our situation since we may assume $e \geq 3$ if $n=5$ by \cite{BC} (resp. $e \geq 8$ if $n=6$ by \cite{HS} Cor. 6.2). So we may assume $d(Z)\cdot (e+n+1) \leq (s-1)[(s-1)^2-d(Z)]\: (*)$. Now if $e \geq E$, from $(*)$ we get: $d(Z) \leq \frac{(s-1)^3}{E+n+s}\:(+)$.
\end{notations}

\begin{proof}[Proof of Theorem \ref{n=s=5}]
Applying $(+)$ with $n=s=5$ and $E=3$ we get $d(Z) \leq 4$, hence $d(Z)=4$ (\cite{R}). Arguing as in \cite{EF} Lemma 2.6, every irreducible component of $Z_{red}$ appears with multiplicity, so $Z$ is either a multiplicity four structure on a linear space or a double structure on a quadric. In both cases it is a complete intersection: in the first case this follows from \cite{Mano} and in the second one, from the fact that $Z$ is given by the Ferrand construction since $emdim(Z_{red})\leq 4$.
\end{proof}

\begin{proof}[Proof of Theorem \ref{n=s=6}]
Applying $(+)$ with $n=s=6$ and $E=8$, we get $d(Z) \leq 6$. If $d(Z)=6$, $(*)$ implies $e \leq 8$. So $e=8$ and $6=d(Z)=d-6e-6$. It follows that $d=60$ and we conclude with \cite{EF} Theorem 1.1. So $d(Z) \leq 5$, hence (\cite{R}), $d(Z)=5$. Now $(*)$ yields $e \leq 13$. Moreover $5=d(Z)=d-6e-6$ yields $d=6e+11$. If $e \leq 10$, again, we conclude with Theorem 1.1 of \cite{EF}. We are left with the following possibilities: $(d,e)=(77,11),(83,12),(89,13)$. We conclude with \cite{HS} (list on page 216).
\end{proof}

\end{document}